%% file: main.tex
\documentclass{article}
\usepackage{spconf,amsmath,graphicx,amsthm,amssymb,color}
\usepackage{booktabs}
\usepackage{url}
 
\makeatletter
\renewcommand\normalsize{%
 \abovedisplayskip 2.8\p@ \@plus2.8\p@ \@minus2.8\p@
 \belowdisplayskip \abovedisplayskip
 \let\@listi\@listI}
\makeatother
\usepackage[ruled, vlined]{algorithm2e}
\usepackage{bm}
\usepackage{subfiles} 
\usepackage[numbers,sort&compress]{natbib}
\usepackage{enumitem}
\usepackage{paralist}
\usepackage{hyperref}
\hypersetup{colorlinks={true},linkcolor={magenta},citecolor={blue}}

\newcommand{\blue}{}

\makeatletter
\renewcommand\normalsize{%
\@setfontsize\normalsize\@xpt\@xiipt
\abovedisplayskip 8\p@ \@plus3\p@ \@minus5\p@
\abovedisplayshortskip \z@ \@plus3\p@
\belowdisplayshortskip 6\p@ \@plus3\p@ \@minus3\p@
\belowdisplayskip \abovedisplayskip
\let\@listi\@listI}
\makeatother
\title{An Online Algorithm for Chance constrained Resource Allocation}
\name{Yuwei Chen$^{*1}$\thanks{$^*$ Denotes equal contribution. The corresponding author is Yuwei Chen
(chen\_yw16@163.com).}
, Zengde Deng$^{*1}$, Yinzhi Zhou$^1$, Zaiyi Chen$^1$, Yujie Chen$^1$, and Haoyuan Hu$^1$}
\address{$^1$Cainiao Network, Hangzhou, China. }
\begin{document}
\ninept
\maketitle
\begin{abstract}
This paper studies the online stochastic resource allocation problem~(RAP) with chance constraints. 
The online RAP is a 0-1 integer linear programming problem where the resource consumption coefficients are revealed column by column along with the corresponding revenue coefficients. 
When a column is revealed, the corresponding decision variables are determined instantaneously without future information. Moreover, in online applications, the resource consumption coefficients are often obtained by prediction. To model their uncertainties, we take the chance constraints into the consideration.
To the best of our knowledge, this is the first time chance constraints are introduced in the online RAP problem. Assuming that the uncertain variables have known Gaussian distributions, the stochastic RAP can be transformed into a deterministic but nonlinear problem with integer second-order cone constraints. 
Next, we linearize this nonlinear problem and analyze the performance of vanilla online primal-dual algorithm for solving the linearized stochastic RAP. Under mild technical assumptions, the optimality gap and constraint violation are both on the order of $\sqrt{n}$.
Then, to further improve the performance of the algorithm, several modified online primal-dual algorithms with heuristic corrections are proposed. 
Finally, extensive numerical experiments \blue{on both synthetic and real data} demonstrate the applicability and effectiveness of our methods.
\end{abstract}
\vspace{-0.3em}
\section{Introduction}
\vspace{-0.3em}
\label{sec:intro}
\input{sections/introduction}

\vspace{-0.3em}
\section{Model Description}
\vspace{-0.4em}
\label{model}
\input{sections/model}

\section{Solution Algorithms}
\label{algorithm}
\input{sections/algorithm}

\section{Numerical Experiments}
\label{exper}
\input{sections/experiments}

\section{Conclusion}
\label{conclusion}
\input{sections/conclusion}

%


\vfill\pagebreak

\bibliographystyle{IEEEbib}
\bibliography{refs}

\end{document}

%% file: sections/introduction.tex
The resource allocation problem~(RAP) \cite{asadpour2020online} is to find the best allocation of a fixed amount of resources to various activates, in order to maximize the total revenue.
The online RAP has a wide range of applications such as signal processing \cite{5230846}, computer resource allocation \cite{kurose1989microeconomic} and
portfolio selection \cite{ida2003portfolio}. 
\blue{This paper studies a multi-dimensional online RAP with uncertainty.
There are $m$ resources and $k$ resource consumption schemes for each request. The request for the resources arrives one by one. When the $i$-th request is revealed, one or none of $k$ resource consumption schemes is chosen to satisfy this request. 
If $l$-th resource consumption scheme is chosen, the revenue and the consumption of the $j$-th resource are $c_{tl}$ and $a_{tjl}$ respectively. 
The decision is irrevocable and has to be decided immediately according to the historical information $\{(\bm c_\tau, \bm A_\tau)\}_{\tau=1}^t$, without future information. 
Our aim is to maximize the total revenue with limited resource capacities, given the total number $n$ of incoming requests and considering the uncertainty of $a_{tjl}$.}

The deterministic RAP can be modeled as a 0-1 integer linear programming~(ILP) problem. 
Many recent papers have studied the online ILP problems (see \citep{molinaro2014geometry, li2019online, agrawal2014dynamic, gupta2014experts, chen2015dynamic, gao2021boosting, li2020simple, balseiro2020dual,balseiro2022best} and references therein). \blue{Algorithms in \citep{molinaro2014geometry, li2019online, agrawal2014dynamic, gupta2014experts, chen2015dynamic, gao2021boosting, li2020simple,balseiro2020dual,balseiro2022best} are all dual-based which maintain dual prices in iterations and can achieve near-optimal solutions under mild conditions. 
When a new request arrives, the decision is made immediately based on the dual price vector.
Among these studies, researchers \citep{molinaro2014geometry, li2019online, agrawal2014dynamic, gupta2014experts, chen2015dynamic} construct dual problems by using historical information and solving them to obtain the dual prices. To deal with the disadvantage that solving dual problems may be time-consuming, researchers \citep{gao2021boosting, li2020simple, balseiro2020dual,balseiro2022best} 
propose online primal-dual~(OPD) algorithms that update the dual prices by utilizing the dual mirror descent or projected stochastic subgradient descent without solving optimization problems.
}

However, the optimization models studied in \citep{molinaro2014geometry, li2019online, agrawal2014dynamic,gupta2014experts, chen2015dynamic, gao2021boosting,li2020simple,balseiro2020dual,balseiro2022best} are deterministic and may suffer from poor performance when the resource consumption is uncertain in practice.
In the existing articles that study the uncertain online optimization, the uncertainty is modeled by the worst-case scenario value, expectation, regret, or a linear combination of the above (see \citep{bent2004online, bent2005online, jiang2020online, liu2015averaging, jiang2020online2}).
These modeling methods are mainly aimed at the uncertainty in the objective function, while almost no chance constraint is considered in the existing studies.

Chance constrained programming~(CCP) \citep{charnes1959chance} is a widely used stochastic programming technique \blue{to model the uncertainty in constraints}.
In stochastic programming \citep{ruszczynski2003stochastic,prekopa2013stochastic}, it is assumed that some parameters are uncertain and their distributions are known.
If the uncertain parameters in an active inequality constraint are set to the medians, the probability of this constraint not holding is 50\%. 
To avoid this issue in the online RAP, this paper adopts the chance constraints to model the uncertainty. The chance constraint is the constraint on the uncertain parameters whose holding probability is not lower than the prescribed level. 
The solution methods for CCPs have been studied by \citep{charnes1959chance, jagannathan1974chance, charnes1963deterministic, prekopa1973contributions,li2008chance}. 
If the uncertain parameters have a known multivariate Gaussian distribution, the chance constrained counterparts of linear constraints can be transformed into deterministic second-order cone (SOC) constraints. prekopa2013stochastic

This paper studies the online stochastic RAP considering the uncertainty of resource consumption coefficients. 
The chance constraint are used to model the uncertainty and can be transformed into the SOC constraints equivalently.
The non-linearity and indecomposability of the SOC constraints make the online problem challenging to handle.
The main contributions of this paper are as follows.

\begin{enumerate}
	\item[(1)] To the best of our knowledge, this is the first time chance constraints are introduced in the online RAP. A linearization method is presented to transform the SOC constrained problem into a linear form suitable for the online solution.
	\item[(2)] We theoretically analyze the performance of the vanilla OPD algorithm when it is applied to solve the SOC constrained RAP. Under mild technical assumptions, the expected optimality gap and constraint violation are both $O(\sqrt{n})$.
	\item[(3)] We propose modified versions of the OPD approach by leveraging the structure of the SOC constraints to effectively reduce the probability deviation \blue{in practice}.
	\item[(4)] Massive numerical experiments based on both synthetic and real data are conducted to demonstrate the applicability and effectiveness of the proposed algorithms.
\end{enumerate}


%% file: sections/model.tex
In this section, we first formulate the deterministic model of the RAP. 
Then, a nonlinear chance constrained counterpart of the RAP is established. Finally, the CCP problem is relaxed into an integer linear problem suitable for online solution.

\vspace{-0.5em}
\subsection{Deterministic Problem}
\vspace{-0.5em}
Consider the multi-dimensional RAP with $n$ requests and $m$ resources.
For each request, there are always $k$ resource consumption schemes that can satisfy it. 
When a request is revealed, the decision maker chooses one scheme or none.
Without loss of generality, a deterministic multi-dimensional RAP can be modeled as follows:
\begin{equation}\label{prob:deterministic}
	\setlength{\abovedisplayskip}{3pt}
	\setlength{\belowdisplayskip}{3pt}
	\begin{aligned}
		\max_{\bm{x}}\quad &\sum_{t=1}^n \bm{c}_t^{\top} \bm{x}_t \\
		\rm{s.t.}\quad&\sum_{t=1}^n \bm{a}_{tj}^{\top} \bm{x}_t \le b_j, \forall j = 1,\dots,m\\
		&\bm{1}^{\top} \bm{x}_t \le 1, \bm{x}_t \in \{0,1\}^k, \forall t = 1,\dots,n
	\end{aligned}
\end{equation}
where the revenue coefficient vector $\bm{c}_t \in \mathbb{R}^k$, and the resource consumption vector $\bm{a}_{tj}\in \mathbb{R}^k$. The decision variables are  $(\bm{x}_1,\dots,\bm{x}_n)$. 
$x_{tl} = 1$ means that $t$-th request is satisfied by resource consumption scheme $l$.
$b_j$ is the capacity of resource $j$. $\bm{1}$ denotes all-one vector.
In the online setting of ILP, the input data $(\bm{c}_t,\bm{a}_{t1},\dots,\bm{a}_{tm})$ is revealed one by one and $\bm{x}_t$ is determined instantaneously when $(\bm{c}_t,\bm{a}_{t1},\dots,\bm{a}_{tm})$ is revealed without future information. Moreover, $n$ and $\bm b$ is known and fixed before the first input arrives.

\vspace{-0.6em}
\subsection{CCP Problem}
\vspace{-0.4em}
{In practical, the value of $\bm{a}_{tj}$ can be obtained by prediction which yields the uncertainty.
Consequently, taken the uncertainty of $\bm{a}_{tj}$ into consideration, we formulate the following CCP problem}:
\begin{equation}\label{OriginalCCP}
\setlength{\abovedisplayskip}{3pt} 
\setlength{\belowdisplayskip}{3pt}
	\begin{aligned}
		\max_{\boldsymbol{x}}\quad &\sum\nolimits_{t=1}^n \boldsymbol{c}_t^{\top} \boldsymbol{x}_t \\
		\rm{s.t.}\quad& \mathbb{P}\left(\sum\nolimits_{t=1}^n \boldsymbol{a}_{tj}^{\top} \boldsymbol{x}_t \le b_j\right)\ge \eta_j,\forall j = 1,\dots,m\\
		&\boldsymbol{1}^{\top} \boldsymbol{x}_t \le 1, \bm{x}_t \in \{0,1\}^k, \forall t = 1,\dots,n
	\end{aligned}
\end{equation}
where $\mathbb{P}$ means probability and $\eta_j$ is the given confidence level.
Assume that the true value of $\bm{a}_{tj}$ belongs to a known Gaussian distribution with mean $\bar{\boldsymbol{a}}_{tj}$ and covariance matrix $\bm{K}_{tj}$ \cite{HOSSEINZADEHLOTFI20121783}, and then problem \eqref{OriginalCCP} is equivalent to the following deterministic problem \cite{boyd2004convex}:
\begin{equation}\label{NLCCP}
\setlength{\abovedisplayskip}{3pt} 
\setlength{\belowdisplayskip}{3pt}
	\begin{aligned}
		\max_{\boldsymbol{x}}\quad &\sum\nolimits_{t=1}^n \boldsymbol{c}_t^\top \boldsymbol{x}_t \\
		\rm{s.t.}\quad  &\sum\nolimits_{t=1}^n \bar{\boldsymbol{a}}_{tj}^\top \boldsymbol{x}_t + \Phi^{-1}(\eta_j)\sqrt{\sum\nolimits_{t=1}^n \boldsymbol{x}_t^\top \bm{K}_{tj} \boldsymbol{x}_t}\\
		& \le b_j,\forall j = 1,\dots,m\\
		&\boldsymbol{1}^{\top} \boldsymbol{x}_t \le 1, \bm{x}_t \in \{0,1\}^k, \forall t = 1,\dots,n
	\end{aligned}
\end{equation}
where $\Phi(\cdot)$ represents the cumulative distribution function of the standard Gaussian
distribution. 
{Moreover, when $\bm{a}_{tj}$ follows a distribution with a finite support, which need not to be Gaussian distribution, problem (\ref{OriginalCCP}) can also be translated into the same form in (\ref{NLCCP}) according to the previous work~\cite{cohen2019overcommitment}.}
Problem (\ref{NLCCP}) is an integer second-order cone programming~(ISOCP) problem when $\eta_j>50\%,\forall j$. 
The offline ISOCP problems can be solved by commercial solvers such as Gurobi. 
However, in the online setting, it is difficult to solve problem (\ref{NLCCP}) due to its non-linearity:  $\bm{x}_t$ with different subscripts {\it t} are coupled with each other in $({\sum_{t=1}^n \boldsymbol{x}_t^\top \bm{K}_{tj} \boldsymbol{x}_t})^{1/2}$. In the online setting of CCP, the input data is $(\bm{c}_t,\bar{\bm{a}}_{t1},\dots,\bar{\bm{a}}_{tm},\bm{K}_{t1},\dots,\bm{K}_{tm})$.

\vspace{-0.6em}
\subsection{Relaxed Linear Problem} 
\newtheorem{prop}{Proposition}[]
\begin{prop}\label{prop:1}
For all $t$ and $j$, the following equation holds.
\begin{equation*}
\setlength{\abovedisplayskip}{2pt} 
\setlength{\belowdisplayskip}{2pt}
	\sqrt{\boldsymbol{x}_t^\top \bm{K}_{tj} \boldsymbol{x}_t} = \boldsymbol{\bm{\gamma}}^\top_{tj} \boldsymbol{x}_t,\forall \boldsymbol{x}_t\in\{\boldsymbol{x} \in \{0,1\}^k|\boldsymbol{1}^\top\boldsymbol{x} \le 1\},
\end{equation*}
where $\bm{\gamma}_{tj}$ is formed by concatenating the square roots of the diagonal elements of the matrix $\bm{K}_{tj}$.
\end{prop}
\vspace{-0.3em}
{To address the non-decomposable issue raised by the non-linearity of $\sqrt{\sum_{t=1}^n \boldsymbol{x}_t^\top \bm{K}_{tj} \boldsymbol{x}_t}$, we linearize this term to decouple different $\bm{x}_t$.} 
Specifically, according to Cauchy-Schwarz inequality
$
	\sqrt{n\sum_{t=1}^n\boldsymbol{x}_t^\top \bm{K}_{tj} \boldsymbol{x}_t} \ge \sum_{t=1}^n \sqrt{\boldsymbol{x}_t^\top \bm{K}_{tj} \boldsymbol{x}_t}
$
and Proposition \ref{prop:1}, the nonlinear problem (\ref{NLCCP}) can be approximated by
\begin{equation} \label{RelaxedCCP}
\setlength{\abovedisplayskip}{3pt} 
\setlength{\belowdisplayskip}{3pt}
	\begin{aligned}
		\max_{\bm{x}}\quad &\sum\nolimits_{t = 1}^n \bm{c}_t^\top \bm{x}_t \\
		\rm{s.t.}\quad&\sum\nolimits_{t = 1}^n \big(\bar{\bm{a}}_{tj}^\top + {\Phi^{-1}(\eta_j)}{} \bm{\gamma}^\top_{tj}/\sqrt{n}  \big) \bm{x}_t \le b_j,\\
		&\forall j = 1,\dots,m\\
		&\boldsymbol{1}^{\top} \boldsymbol{x}_t \le 1, \bm{x}_t \in \{0,1\}^k, \forall t = 1,\dots,n.
	\end{aligned}
\end{equation}
Problem (\ref{RelaxedCCP}) is linear and can be solved in the online setting. 
The online algorithm for solving this relaxed problem is the basis of our algorithm for solving the CCP problem (\ref{NLCCP}) and we will detail it in the following section.


%% file: sections/algorithm.tex
\blue{In this section, we introduce several online primal-dual methods to handle the online SOC constrained problem \eqref{NLCCP}. Firstly, we revisit the state-of-the-art OPD algorithm for solving the relaxed problem \eqref{RelaxedCCP}. Then, some heuristic correction methods based on the structure of \eqref{NLCCP} are proposed to improve the practical performance.}

\vspace{-0.6em}
\subsection{OPD Algorithm for online ILP}
\blue{Recall that \eqref{RelaxedCCP} is an ILP problem} and Li et al.~\cite{li2020simple} have proposed an effective OPD algorithm \blue{to solve the online ILP problem}. 
\blue{For simplicity, denote $\tilde{\bm{a}}_{tj} = \bar{\bm{a}}_{tj} + {\Phi^{-1}(\eta_j)}\bm{\gamma}_{tj}{/\sqrt{n}}$ and we present the OPD method as shown in Algorithm \ref{Algo 1}.}

\begin{algorithm}[tb]\label{Algo 1}
\caption{OPD Algorithm for ILP}
\LinesNumbered
\KwIn{$\bm{d} = \bm{b}/n$}
\KwOut{$\bm{x}=(\bm{x}_1,...,\bm{x}_n)$}
{\bf Initialize:} $\bm{p}_1=\bm{0}$\\
\For{$t=1,...,n$}
{
	Set $v_t = \max_{l=1,\dots,k} \  (\bm{c}_t^\top - \bm{p}_t^\top \tilde{\bm{A}}_{t})\bm{e}_l$\\
	\eIf{$v_t > 0$}
	{
		Pick an index $l_t$ randomly from
 \begin{equation*}
	\setlength{\abovedisplayskip}{0pt}
	\setlength{\belowdisplayskip}{0pt}
 	 \big\{l:v_t=(\bm{c}_t^\top - \bm{p}_t^\top \tilde{\bm{A}}_{t})\bm{e}_l\big\}
 \end{equation*}\\	
  		Set $\bm{x}_t = \bm{e}_{l_t}$
	}{Set $\bm{x}_t = \bm{0}$}
	Compute 
	{
	$\bm{p}_{t+1} = \max\Big\{\bm{p}_{t} + \frac{1}{\sqrt{n}}\big(\tilde{\bm{A}}_{t}\bm{x}_t-\bm{d}\big), \mathbf{0}\Big\}$
	}
}
\end{algorithm}

In Algorithm \ref{Algo 1}, denote $\tilde{\bm{A}}_t = (\tilde{\bm{a}}_{t1}^\top,\dots,\tilde{\bm{a}}_{tm}^\top)^\top$ and $\bm{b} = (b_1,\dots,b_m)^\top$. 
Algorithm \ref{Algo 1} is dual-based which maintains a dual vector $\bm{p}_t$. 
In each iteration, new $\bm{c}_t$ and $\tilde{\bm{A}}_{t}$ are revealed. 
Then, $\bm{x}_t$ is determined by choosing $l$ that maximizes $(\bm{c}_t^\top - \bm{p}_t^\top \tilde{\bm{A}}_{t})\bm{e}_l$, where $\bm{e}_l$ is the unit vector with all components equal to 0 except the $l$-th, which is 1. 
{After determining $\bm{x}_t$, $\bm{p}_t$ is updated by a projected stochastic subgradient descent method where $(\bm{d} - \tilde{\bm{A}}_{t}\bm{x}_t)$ is the subgradient corresponding to $\bm{p}_t$.}

The following Theorem \ref{thm: regret and violation ISOCP} states that  Algorithm \ref{Algo 1} achieves $O(\sqrt{n})$ regret and constraint violation compared to the optimal solution of the ISOCP problem \eqref{NLCCP}.
The detailed proof of Theorem \ref{thm: regret and violation ISOCP} as well as Proposition \ref{prop:1} is presented in the full-length version \cite{chen2022online}.

\newtheorem{thm}{\bf Theorem}[]
\begin{thm}\label{thm: regret and violation ISOCP}
Assume coefficient sets $\{c_{tj}, \bm{\bar{a}}_{tj}, \bm K_{tj}\}$s are bounded and sampled i.i.d.~from an unknown distribution, and the upper and lower bounds of $\bm{b}/n$ are finite and positive. Then, the expected regret and constraint violation of Algorithm \ref{Algo 1} compared to the optimal solution of the ISOCP problem \eqref{NLCCP} are on the order of $\sqrt{n}$, i.e.,
\begin{equation}\label{eq:regret}
\setlength{\abovedisplayskip}{3pt} 
\setlength{\belowdisplayskip}{3pt}
   {\mathbb E_{\{c_{tj}, \bm{\bar{a}}_{tj}, \bm K_{tj}\}_{t=1}^n}} \bigg[\hat{R}_n^{ISOCP}-\sum_{t = 1}^n \bm{c}_t^\top \bm{x}_t\bigg] \le O(\sqrt{n})
\end{equation}
\begin{equation} 
    {\mathbb E_{\{c_{tj}, \bm{\bar{a}}_{tj}, \bm K_{tj}\}_{t=1}^n}} \left[\left\|\left(\bm g\left({{\bm{x}}}\right)-\bm{b}\right)^+\right\|_2\right] \le O(\sqrt{n})
\end{equation}
where $\hat{R}_n^{ISOCP}$ is the optimal objective value of \eqref{NLCCP}, $\bm{x}=(\bm{x}_1,...,\bm{x}_n)$ is the output of Algorithm \ref{Algo 1}, $(\cdot)^+$ is the positive part function, and $g(\bm x)$ is the left-hand side of the SOC constraints. 
\end{thm}

\vspace{-1em}
\subsection{Modified OPD Algorithm for online CCP}
\vspace{-0.4em}
Although Algorithm \ref{Algo 1} has been able to obtain a near-optimal solution of the ISOCP problem \eqref{NLCCP} according to Theorem \ref{thm: regret and violation ISOCP}, its practical performance can be further improved \blue{by narrowing the gap between the solutions generated by Algorithm~\ref{Algo 1} and the offline ISOCP \eqref{NLCCP}}. 
\blue{To be specific,} this gap mainly comes from the following two points:
\begin{enumerate}
\renewcommand{\labelenumi}{(\alph{enumi})}
	\item The error between the offline ILP problem \eqref{RelaxedCCP} and the offline ISOCP problem \eqref{NLCCP}.
	\item The error between the online solution and offline solution of the ILP problem \eqref{RelaxedCCP}.
\end{enumerate}

\begin{algorithm}[t]\label{Algo 2}
\caption{Modified OPD Algorithm for CCP}
\LinesNumbered
\KwIn{$\bm{d} = \bm{b}/n$}
\KwOut{$\bm{x}=(\bm{x}_1,...,\bm{x}_n)$}
{\bf Initialize:} $\bm{p}_1=\bm{0}$, $\bm{d}_1 = \bm{d}$\\
\For{$t=1,...,n$}
{
	Compute $\bm{\beta}_{t}$ via equation (\ref{beta})\\
	Set $v_t = \max_{l=1,\dots,k} \  (\bm{c}_t^\top - \bm{p}_t^\top \hat{\bm{A}}_{t}(\bm{\beta}_t))\bm{e}_l$\\
	\eIf{$v_t > 0$}
	{
		Pick an index $l_t$ randomly from
 \begin{equation*}
	\setlength{\abovedisplayskip}{0pt}
	\setlength{\belowdisplayskip}{0pt}
	\left\{l:v_t=(\bm{c}_t^\top - \bm{p}_t^\top \hat{\bm{A}}_{t}(\bm{\beta}_t))\bm{e}_l\right\}
	\end{equation*}\\	
 		Set $\bm{x}_t = \bm{e}_{l_t}$
	}{Set $\bm{x}_t = \bm{0}$}
	Compute $\bm{d}_{t}$ via equation \eqref{d_t}\\
	Compute 
	{
	 \begin{equation*}
	\setlength{\abovedisplayskip}{0pt}
	\setlength{\belowdisplayskip}{0pt}
	\ \bm{p}_{t+1} = \max\Big\{\bm{p}_{t} + \frac{1}{\sqrt{n}}\left(\hat{\bm{A}}_{t}(\bm{\beta}_t)\bm{x}_t-\bm{d}_t\right),\mathbf{0}\Big\}
		\end{equation*}
	}
}
\end{algorithm}

To address these issues, we propose the modified OPD Algorithm~\ref{Algo 2} for solving the online CCP problem~(\ref{NLCCP}).
In Algorithm \ref{Algo 2}, several heuristic corrections are applied to correct the above-mentioned errors. 

First is to correct the error (a). \blue{For the $j$-th constraint,} we introduce scale factors 
\begin{equation}\label{beta}
\!\!\!\!
\beta_{tj}=
\begin{cases}
	1, &t = 1 \text{ or } \sum_{i=1}^{t-1} \bm\gamma_{ij}^\top \bm x_i = 0\\
	\frac{\sqrt{t-1}\sqrt{\sum_{i=1}^{t-1} \bm{x}_i^\top \bm{K}_{ij} \bm{x}_i}}
		{\sum_{i=1}^{t-1} \bm{\gamma}_{ij}^\top \bm{x}_i}, &t > 1\text{ and } \sum_{i=1}^{t-1} \bm\gamma_{ij}^\top \bm x_i > 0
\end{cases}
\end{equation}
to reduce the gap between $\sqrt{\sum_{t=1}^n \boldsymbol{x}_t^\top \bm{K}_{tj} \boldsymbol{x}_t}$ and $\sum_{t=1}^n {\bm{\gamma}_{tj}^\top \boldsymbol{x}_t}$ $/{\sqrt{n}}$. 
Next, we define
$
\hat{\bm{a}}_{t1}(\beta_{tj}) = \bar{\bm{a}}_{tj} + \beta_{tj}{\Phi^{-1}(\eta_j)}{}\bm{\gamma}_{tj}/\sqrt{n}$
and
$
	\hat{\bm{A}}_t(\bm{\beta_t})$ = $(\hat{\bm{a}}_{t1}^\top(\beta_{t1}),\dots,\hat{\bm{a}}_{tm}^\top(\beta_{tm}))^\top.
	$
In Algorithm \ref{Algo 2}, $\hat{\bm{A}}_t(\bm{\beta_t})$ is used in place of $\tilde{\bm{A}}_t$. 
That is, we use $\sum_{t=1}^n \beta_{tj}\bm{\gamma}_{tj}^\top \bm{x}_t/\sqrt{n}$ to approximate $\sqrt{\sum_{t=1}^n \bm{x}_t^\top \bm{K}_{tj} \bm{x}_t}$. At time $t$, $\bm{\beta}_{t}$ is calculated \blue{according to \eqref{beta}} which is based on the historical decisions and will be used in the next iteration for correction.
It is worth noting that $\bm{\beta}_{t} = (\beta_{t1},\dots,\beta_{tm})^\top$ is calculated in each round and can be computed incrementally with low computational cost.

In the numerical experiments section \ref{exper}, it is illustrated that Algorithm \ref{Algo 2} has better performance than Algorithm \ref{Algo 1} in terms of the constraint violation. 
\blue{An intuitive explanation is that Algorithm \ref{Algo 2} is more inclined to reject the orders with high uncertainty of resource consumption (i.e., $\bm{K}_{tj}$) than Algorithm \ref{Algo 1} because $\bm \beta_t \ge \bm 1$.}

Next is to correct the error (b). The error (b) consists of two parts, the optimality gap and constraint violation.
The constraint violation will cause the probability
$
	\mathbb{P}\left(\sum\nolimits_{t=1}^n \boldsymbol{a}_{tj}^{\top} \boldsymbol{x}_t \le b_j \right)
$
to deviate from the target value $\eta_j$. 
It is almost impossible to reduce the optimality gap and probability deviation simultaneously. 
Compared with the optimality gap, the CCP problems have a lower tolerance for the probability deviation. 
In order to reduce the probability deviation, we propose the following method to dynamically adjust the right-hand-side capacity $\bm{d}$ in each iteration:
\begin{equation}\label{d_t}
\setlength{\abovedisplayskip}{2pt} 
\setlength{\belowdisplayskip}{2pt}
\begin{aligned}
		d_{tj} = & \frac{1}{n-t}\bigg(b_j - \Phi^{-1}(\eta_j)\sqrt{\frac{t}{n}\sum\nolimits_{i=1}^t \boldsymbol{x}_i^\top \bm{K}_{ij} \boldsymbol{x}_i} \\
		& -\sum\nolimits_{i=1}^t \bar{\boldsymbol{a}}_{ij}^\top \boldsymbol{x}_i \bigg),\forall j = 1,\dots,m.
\end{aligned}
\end{equation}

\blue{The intuition behind the correction formula \eqref{d_t} is given as follows:} 
if too many resources are spent in the early rounds, the average remaining resources $\bm d$ will diminish. Then Algorithm \ref{Algo 2} will raise the dual price and be more likely to reject an order with high resource consumption as a result. On the other hand, if a large number of orders with high resource consumption are rejected at the start, resulting in an excess of remaining resources, Algorithm \ref{Algo 2} will decrease the dual price in order to accept more orders in the future. This \blue{correction strategy} makes Algorithm \ref{Algo 2} perform better than Algorithms \ref{Algo 1} in numerical experiments.

%

%% file: sections/experiments.tex
In this section, we compare the performance of Algorithm \ref{Algo 1}, Algorithm \ref{Algo 2}, Algorithm \ref{Algo 2} without correction (\ref{beta}) and Algorithm \ref{Algo 2} without correction (\ref{d_t}) in terms of optimality gap and probability deviation. 
These algorithms are implemented on two different models, with details given in Table \ref{tab:model}. 
Table \ref{tab:model} lists the distributions from which the elements in $\bm{c}_{tj}$, $\bm{\bar{a}}_{tj}$ or $\bm{K}_{tj}$ are i.i.d.~sampled in two synthetic-data experiments. $X\sim f(\chi^2(v))$ denotes $X = f(Y)$ and $Y\sim~\chi^2(v)$.
\begin{table}[h]
	\renewcommand\arraystretch{1.35}
	\centering
	\vspace{-0.5cm}
	\caption{Models used in the experiments.}
	\label{tab:model}
	{
		\setlength{\tabcolsep}{4pt}{
			\begin{tabular}{c c c c c c c}
				\toprule[1.5pt]
				Experiment & $\bm{c}_{tj}$ & $\bm{\bar{a}}_{tj}$ & $\bm{K}_{tj}$ & $\bm{d}$\\
				\midrule[0.75pt]
				I & U[0, 1] & U[0, 4] & $(\text{U}[0, 1])^2$ & 1\\
				II & $\chi^2(3)$& $\frac{2}{3}\chi^2(4)$& $(\frac{2}{3}\chi^2(2))^2$ & 1\\
				\bottomrule[1.5pt]
			\end{tabular}
		}
	}
	\vspace{-0.5cm}
\end{table}

\subsection{Synthetic-data Experiment I (Bounded Setting)}
\vspace{-0.4em}
In the first experiment, we set $k = 5$ and $m = 4$. 
The confidence levels of chance constraints are set to (0.65, 0.75, 0.85, 0.95). 
For each value of $n$, we run 20 simulation trials. 
In each trial, coefficients $\bm{c}_{tj}$, $\bm{\bar{a}}_{tj}$ and $\bm{K}_{tj}$ are resampled.

Fig.\ \ref{fig:expr1_algo_comparison} shows the average optimality gap and probability deviation over all the simulation trials. The probability deviation of the whole problem in one trail is an average:
\begin{equation}\label{eq:prob_violation}
\setlength{\abovedisplayskip}{2pt} 
\setlength{\belowdisplayskip}{2pt}
	\frac{1}{m}\sum_{j=1}^m\bigg(\eta_j-\Phi\bigg(\frac{b_j-\sum_{t=1}^n \bar{\boldsymbol{a}}_{tj}^\top\bm{x}_t}{\sqrt{\sum_{t=1}^n \boldsymbol{x}_t^\top \bm{K}_{tj} \boldsymbol{x}_t}}\bigg)\bigg)^+,
\end{equation}
where $\bm{x}_t$ is the output of the algorithms.
{The calculation formula for the expected optimality gap is (\ref{eq:regret}).
From Fig.~\ref{fig:expr1_algo_comparison}, we observe that the optimality gaps of these algorithms are close and Algorithm \ref{Algo 2} has the smallest probability deviation. 
Fig.~\ref{fig:expr1_algo_comparison} (a) also shows that the optimality gap of Algorithm \ref{Algo 2} is on the order of $\sqrt{n}$.
Fig.~\ref{fig:expr1_detailed_violation} presents the probability deviations of each chance constraint of Algorithm \ref{Algo 1} and \ref{Algo 2}.
Fig.~\ref{fig:expr1_algo_comparison} and \ref{fig:expr1_detailed_violation} both illustrate that the proposed two corrections (\ref{beta}) and (\ref{d_t}) can effectively reduce the probability deviation with \blue{minor negative effects} on the optimality gap.

\begin{figure}[tb]
\begin{minipage}[b]{.48\linewidth}
  \centering
  \centerline{\includegraphics[width=4.0cm]{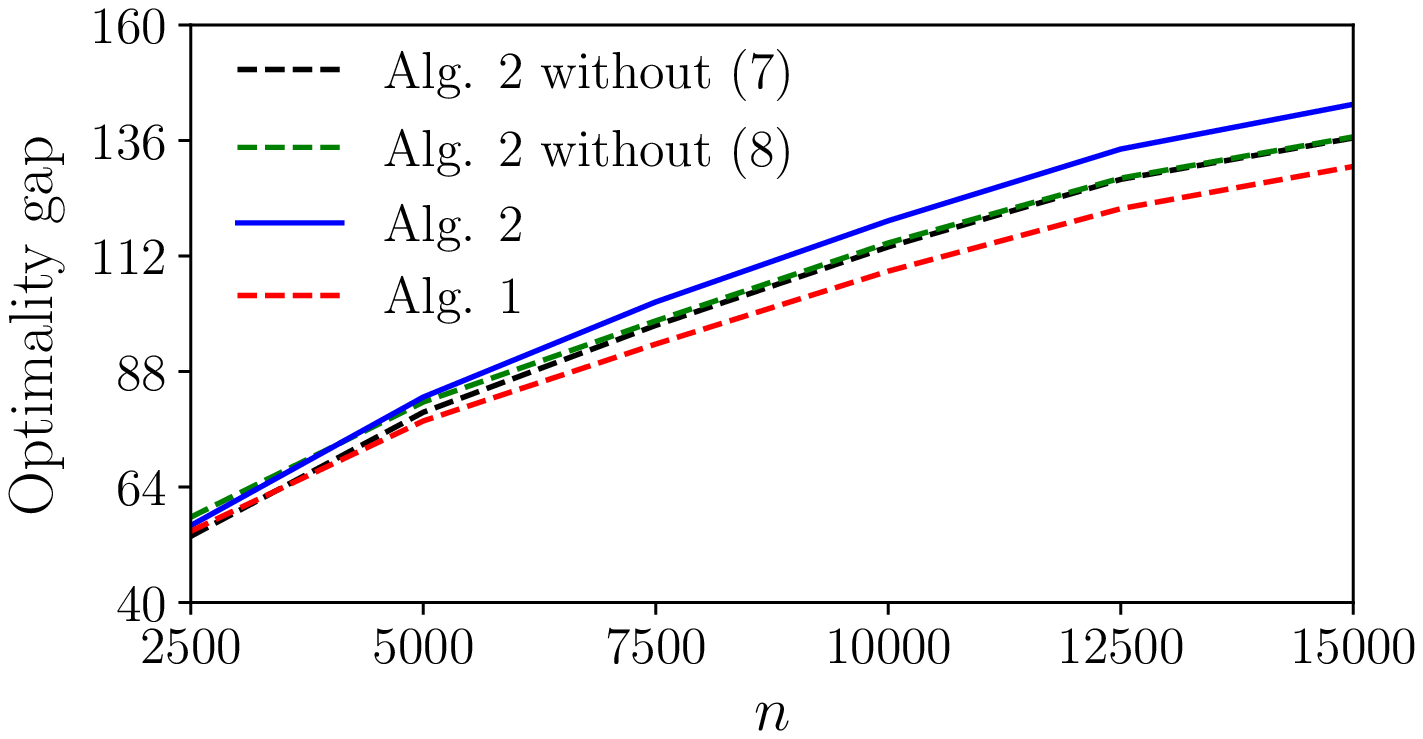}}
     \vspace{-0.01cm}
  \centerline{(a) Optimality gap}\medskip
\end{minipage}
\hfill
\begin{minipage}[b]{0.48\linewidth}
  \centering
  \centerline{\includegraphics[width=4.0cm]{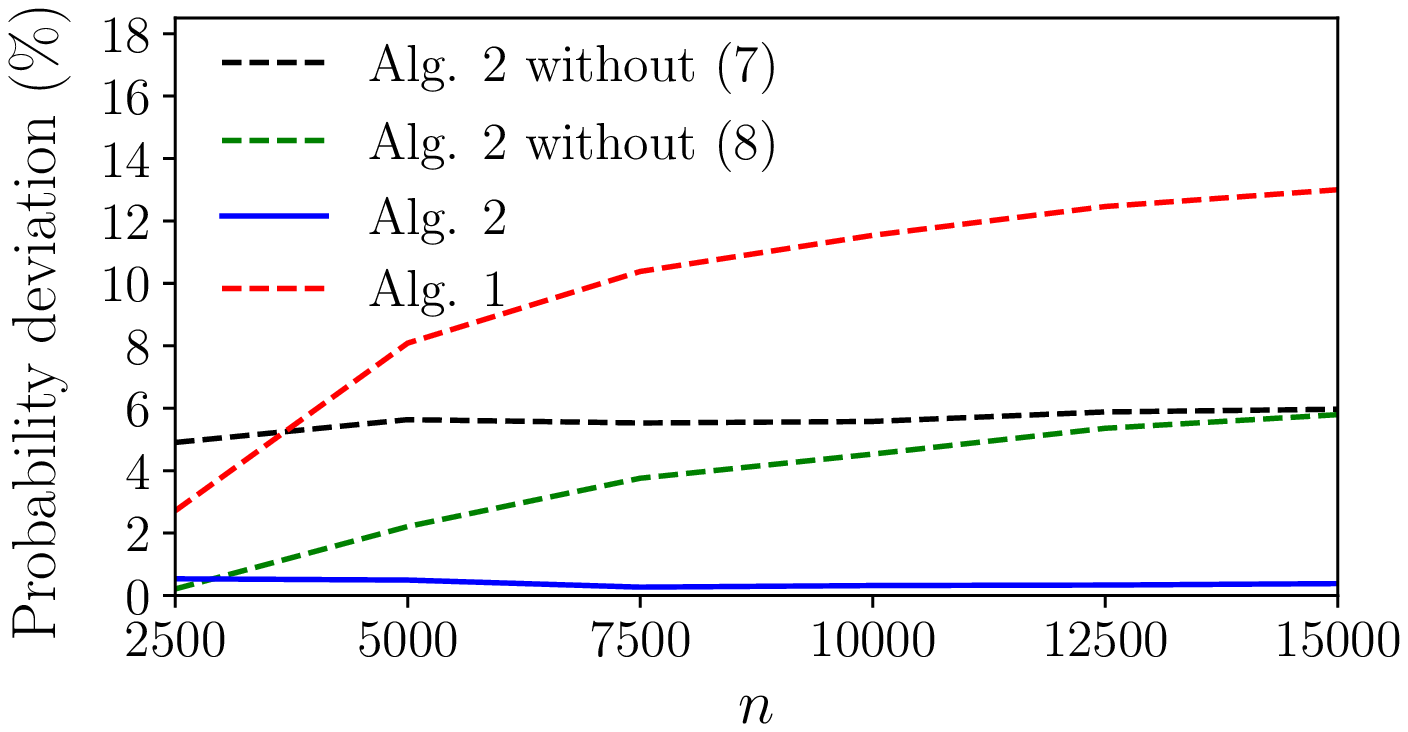}}
     \vspace{-0.01cm}
  \centerline{(b) Probability deviation}\medskip
\end{minipage}
\vspace{-0.5cm}
\caption{Average optimality gap and probability deviation in Experiment I with Uniform i.i.d. input.}
\vspace{-0.2cm}
\label{fig:expr1_algo_comparison}
\end{figure}

\begin{figure}[tb]
\begin{minipage}[b]{.48\linewidth}
  \centering
  \centerline{\includegraphics[width=4.0cm]{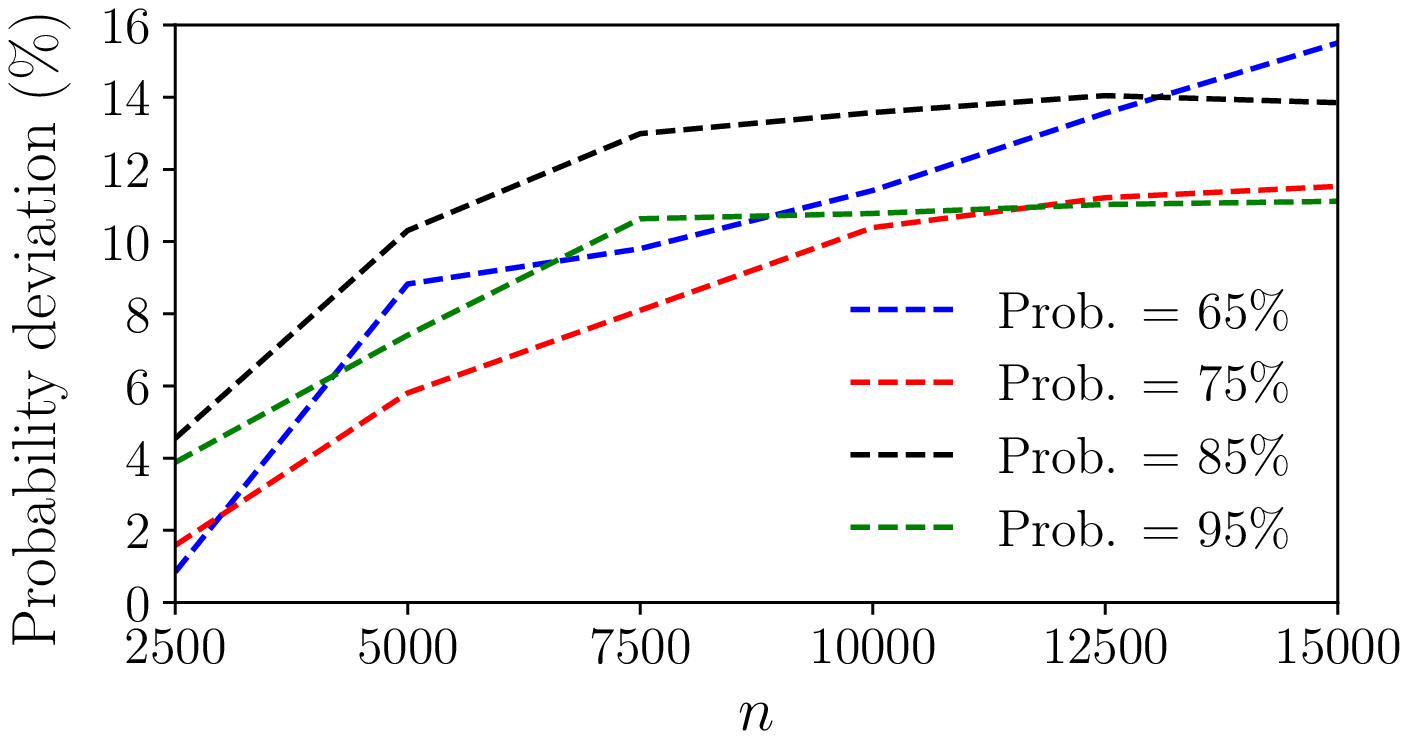}}
     \vspace{-0.01cm}
  \centerline{(a) Algorithm \ref{Algo 1}}\medskip
\end{minipage}
\hfill
\begin{minipage}[b]{0.48\linewidth}
  \centering
  \centerline{\includegraphics[width=4.0cm]{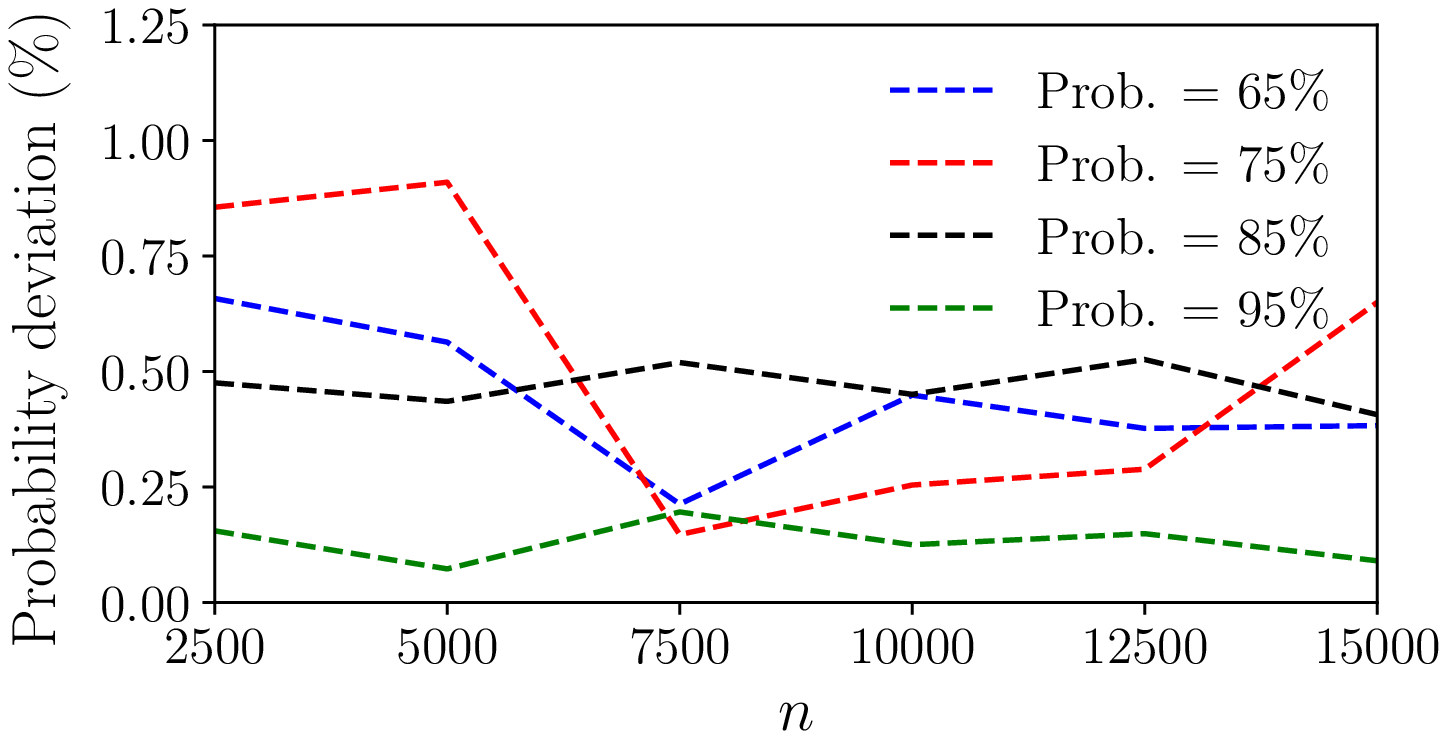}}
    \vspace{-0.05cm}
  \centerline{(b) Algorithm \ref{Algo 2}}\medskip
\end{minipage}
\vspace{-0.5cm}
\caption{Probability deviation of each chance constraint in Experiment I with Uniform i.i.d. input.}
 \vspace{-0.1cm}
\label{fig:expr1_detailed_violation}
\end{figure}

\subsection{Synthetic-data Experiment II (Unbounded Setting)}
\vspace{-0.2em}
In the second experiment, $k$ and $m$ are still set to 5 and 4. The confidence levels are also the same as those in Experiment I. For each value of $n$, we run 20 simulation trials.
In each trial, coefficients $\bm{c}_{tj}$, $\bm{\bar{a}}_{tj}$ and $\bm{K}_{tj}$ are i.i.d.\ sampled from Chi-square distributions which are unbounded.

Fig.~\ref{fig:expr2_algo_comparison} shows the average optimality gap and probability deviation, and Fig.~\ref{fig:expr2_detailed_violation} shows the probability deviations of each chance constraint of Algorithm \ref{Algo 1} and \ref{Algo 2}.
The results of Experiment II are similar to those of Experiment I: 
Algorithm 2 has the smallest probability deviation;
corrections (\ref{beta}) and (\ref{d_t}) can effectively reduce the probability deviation.
Although Algorithm 2 produces slightly larger optimality gap, its optimality gap is still approximately on the order of $\sqrt{n}$. 
In this experiment with unbounded input, Algorithm 2 has obvious advantages:
the probability deviations of the algorithms except Algorithm 2 are larger than 10\%, while the probability deviation of Algorithm \ref{Algo 2} is less than 1\%.

In addition, the competitive ratios of Algorithm 2 in Experiment I and II are not less than 96\%, with details provided in this doc \cite{results2021}.

\begin{figure}[t]
\begin{minipage}[b]{.48\linewidth}
  \centering
  \centerline{\includegraphics[width=4.0cm]{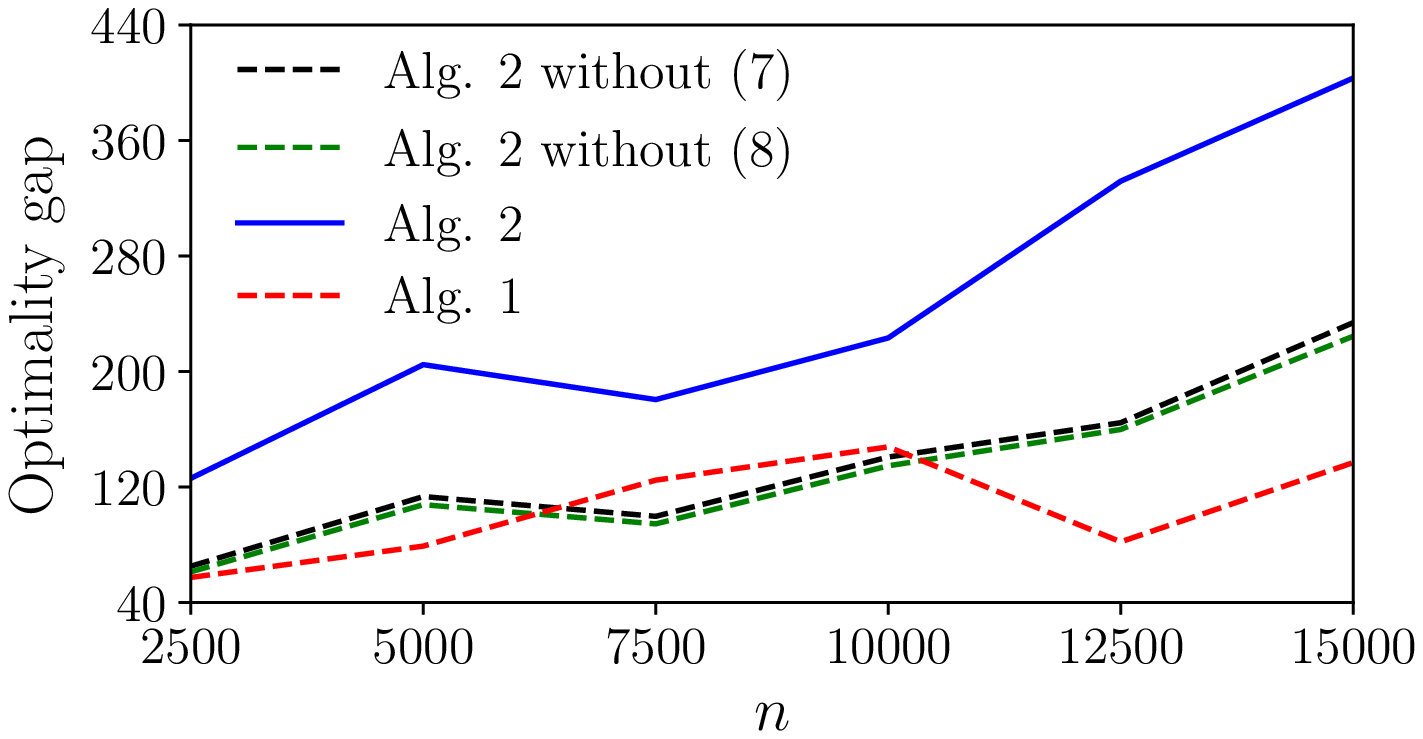}}
       \vspace{-0.01cm}
  \centerline{(a) Optimality gap}\medskip
\end{minipage}
\hfill
\begin{minipage}[b]{0.48\linewidth}
  \centering
  \centerline{\includegraphics[width=4.0cm]{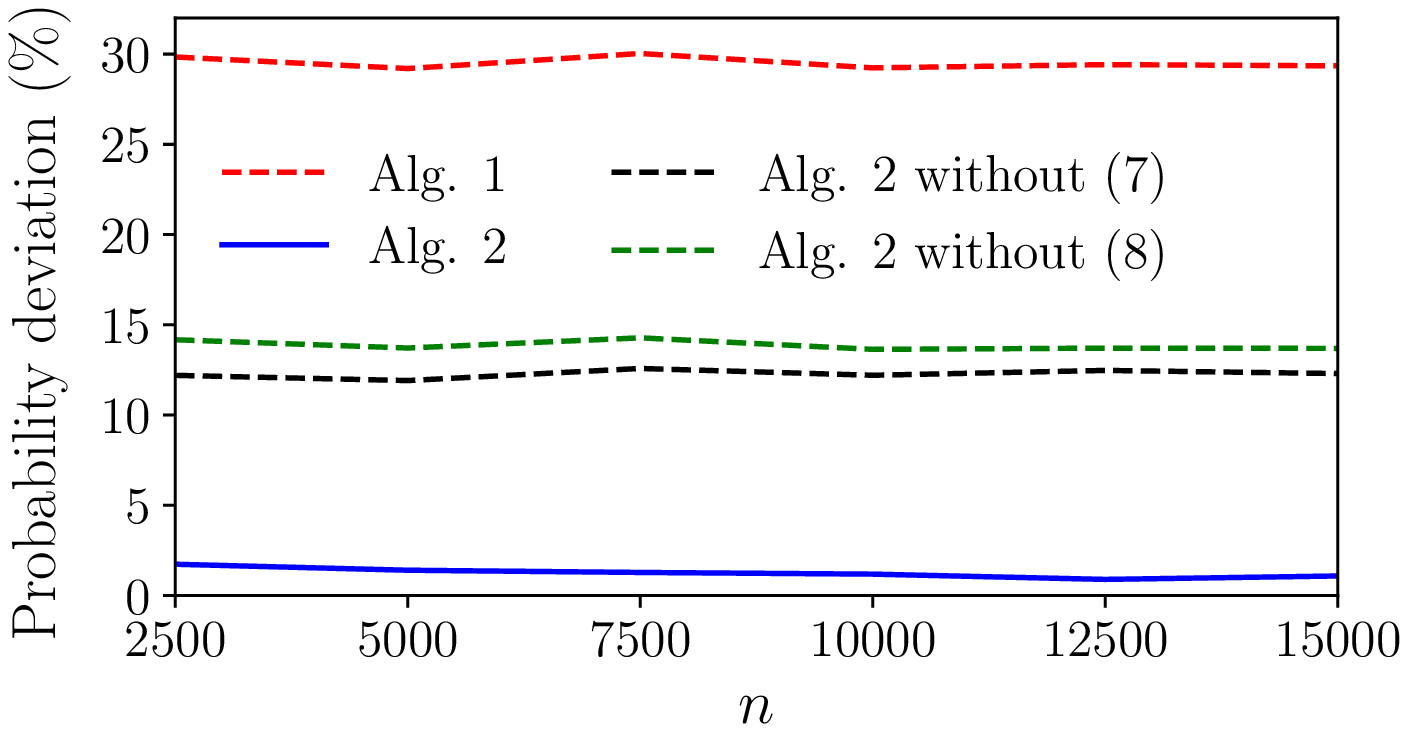}}
       \vspace{-0.01cm}
  \centerline{(b) Probability deviation}\medskip
\end{minipage}
\vspace{-0.5cm}
\caption{Average optimality gap and probability deviation in Experiment II with Chi-square i.i.d. input.}
\vspace{-0.2cm} 
\label{fig:expr2_algo_comparison}
\end{figure}

\begin{figure}[t]
\begin{minipage}[b]{.48\linewidth}
  \centering
  \centerline{\includegraphics[width=4.0cm]{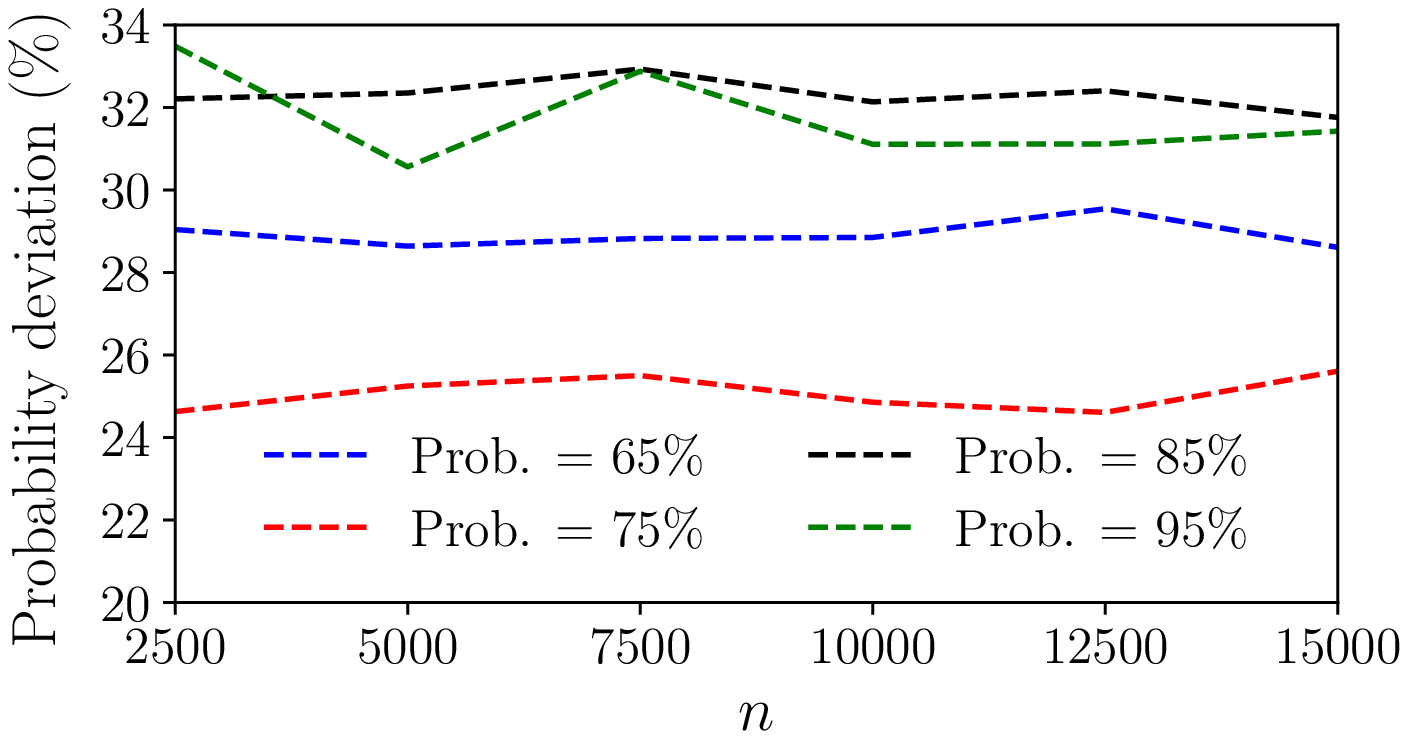}}
       \vspace{-0.01cm}
  \centerline{(a) Algorithm \ref{Algo 1}}\medskip
\end{minipage}
\hfill
\begin{minipage}[b]{0.48\linewidth}
  \centering
  \centerline{\includegraphics[width=4.0cm]{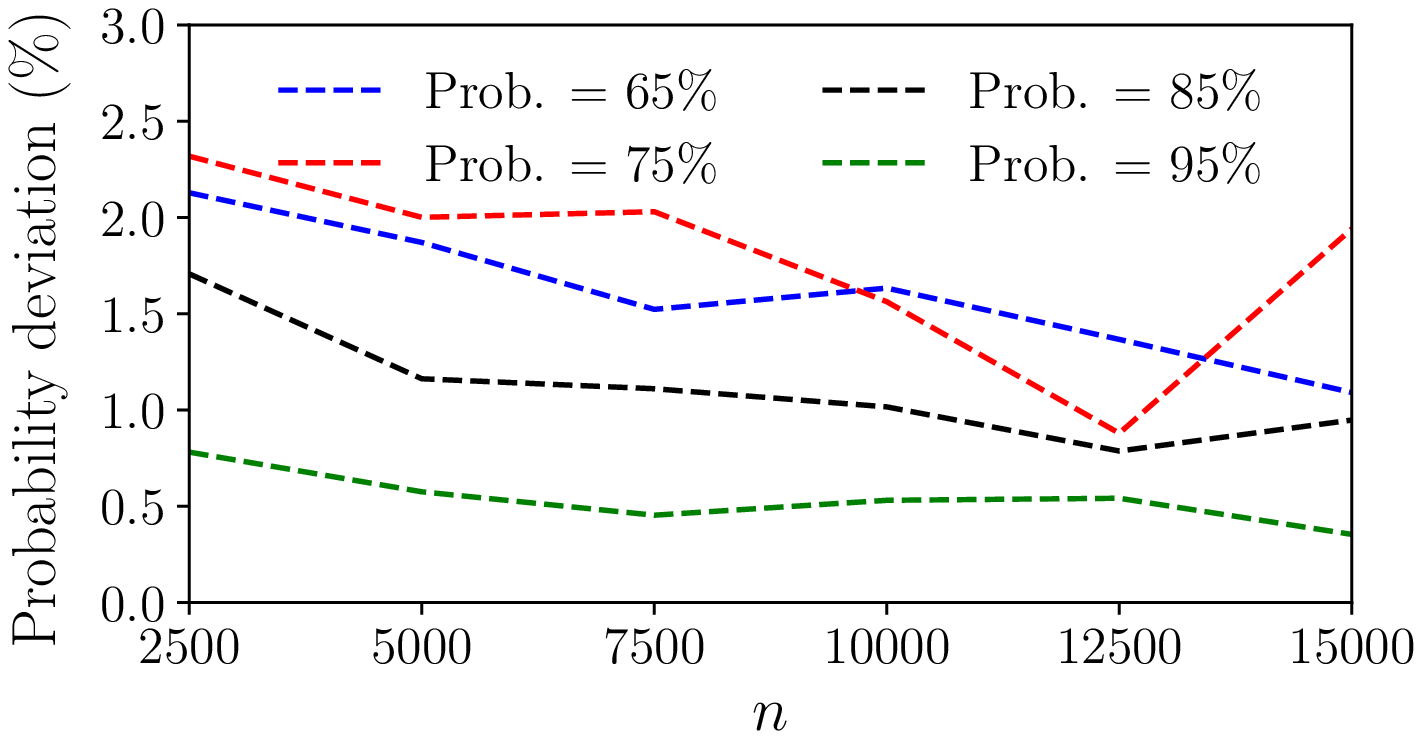}}
       \vspace{-0.01cm}
  \centerline{(b) Algorithm \ref{Algo 2}}\medskip
\end{minipage}
\vspace{-0.5cm}
\caption{Probability deviation of each chance constraint in Experiment II with Chi-square i.i.d. input.}
 \vspace{-0.1cm}
\label{fig:expr2_detailed_violation}
\end{figure}

\subsection{Real-data Experiment}
In the following, we present an engineering application of our method in the task of order fulfillment based on the real data obtained from Cainiao Network \blue{which is a supply chain company}. 
The request (order) is revealed one by one, and the algorithm needs to decide which transportation channel the order will be sent to. 
$x_{tl} = 1$ denotes the order $t$ is sent to channel $l$. 
The objective coefficient is the revenue of each channel and the constraint coefficient is the predicted transportation time of each channel. 
The deterministic problem is to maximize the total revenue while ensuring that the average transportation time is not larger than 15 working days. 
The transportation time is obtained from prediction \blue{which introduces the uncertainty in the constraint coefficients}.
\blue{We adopt the chance constraint with the holding probability of average transportation time $\ge90\%$.}
Due to that each order must be assigned to a channel, the constraint $\bm 1^\top \bm x_t \le 1$ in \eqref{OriginalCCP} is replaced with $\bm 1^\top \bm x_t = 1$ and Algorithms are also slightly modified accordingly: do not judge whether $v_t > 0$ and always set $\bm x_t = \bm e_{l_t}$. Three different transportation channels are considered. 
The simulation results are shown blow.

\begin{figure}[thb]
\vspace{-0.1cm}
\begin{minipage}[b]{.48\linewidth}
  \centering
  \centerline{\includegraphics[width=4.0cm]{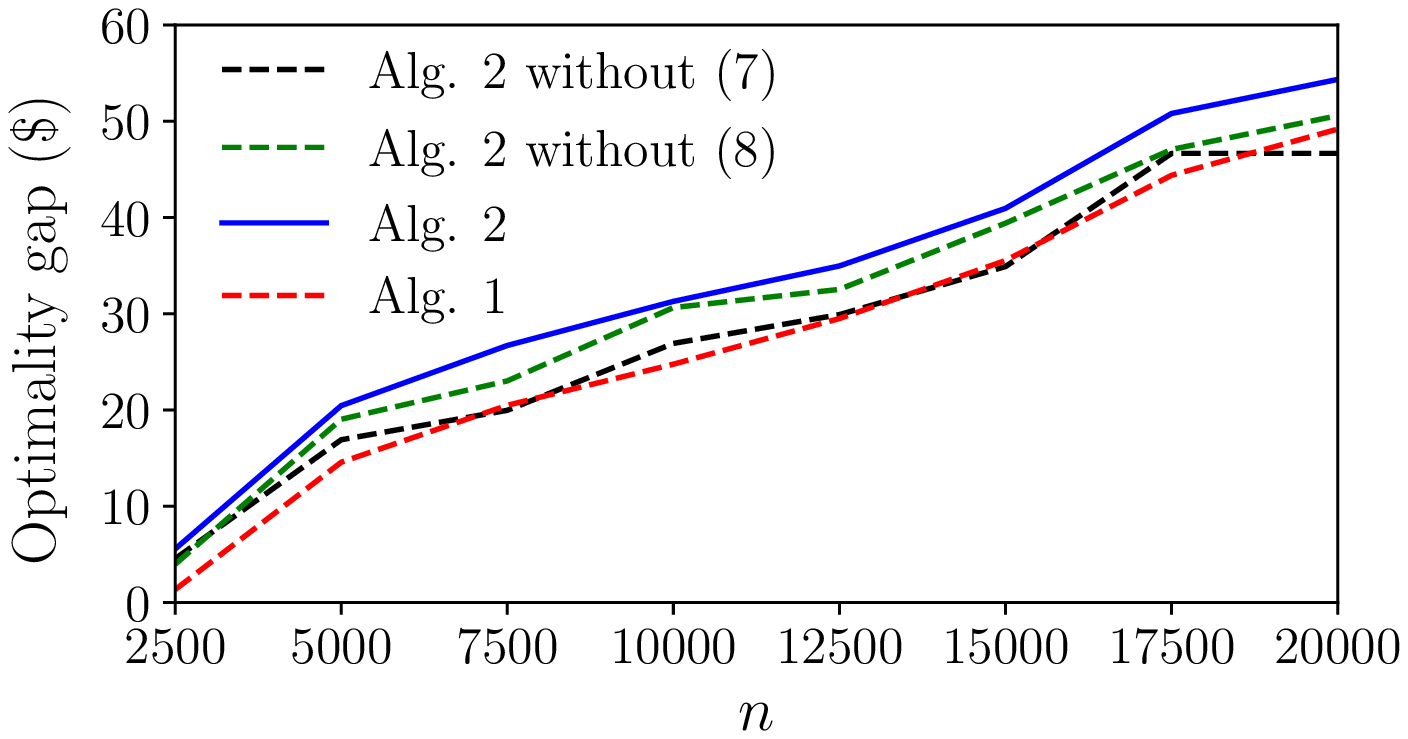}}
       \vspace{-0.01cm}
  \centerline{(a) Optimality gap}\medskip
\end{minipage}
\hfill
\begin{minipage}[b]{0.48\linewidth}
  \centering
  \centerline{\includegraphics[width=4.0cm]{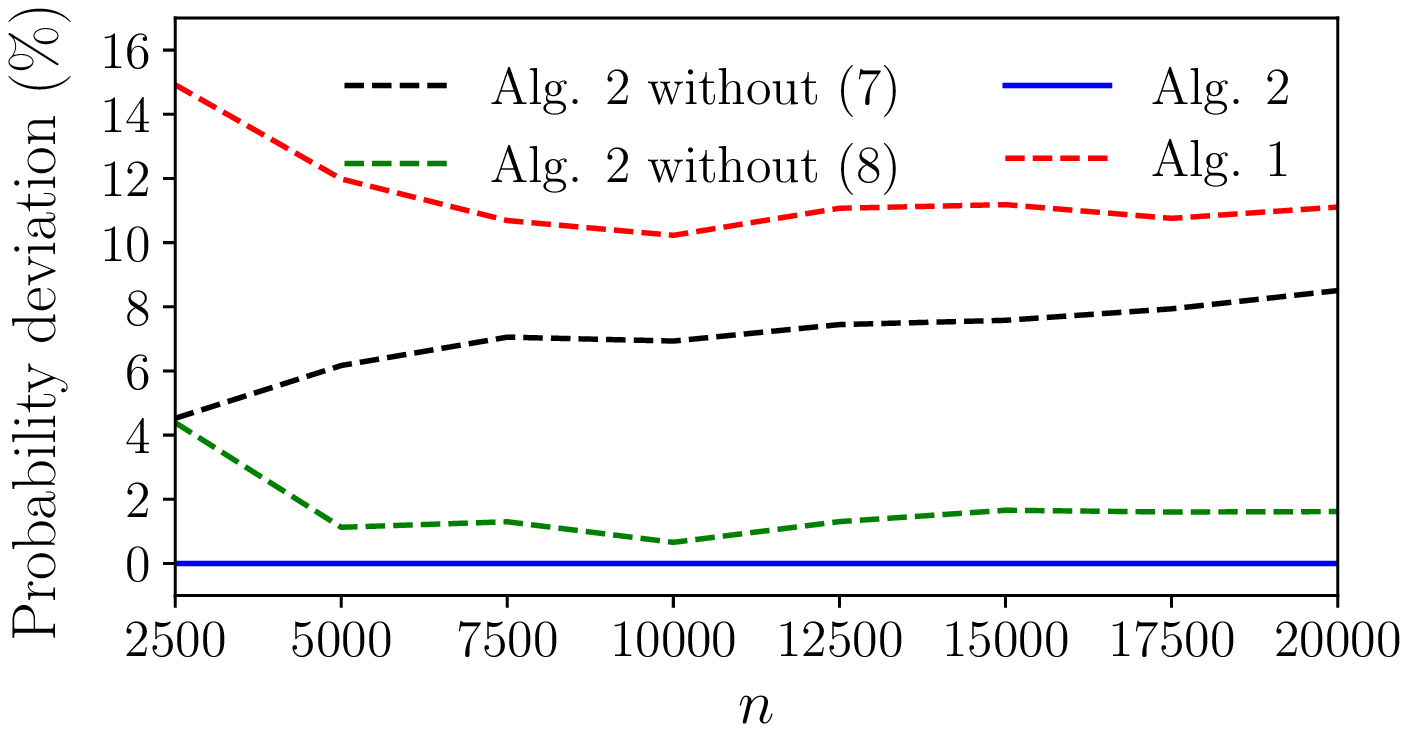}}
       \vspace{-0.01cm}
  \centerline{(b) Probability deviation}\medskip
\end{minipage}
\vspace{-0.5cm}
\caption{Real-data experimental results.}
\vspace{-0.3cm}
    \label{fig:real_expr_algo_comparison}
\end{figure}

The above results also show that Algorithm \ref{Algo 2} with proposed corrections is effective for reducing the value of the probability deviation \blue{while not sacrificing the optimality gap too much}.

%

%% file: sections/conclusion.tex
In this paper, we study the online stochastic RAP with chance constraints. 
First, we present a linearization method that decouples the non-linear term in second-order cone constraints and makes the online solution possible.
Then, we adopt the online primal-dual~(OPD) algorithm for the integer linear programming problem and establish the $O(\sqrt{T})$ regret for both the optimality gap and constraint violation.
Moreover, several heuristic corrections are proposed to further improve the performance of the OPD algorithm.
Extensive numerical experiments on both synthetic and real data verify the effectiveness of our proposed methods.